\documentclass[11pt,letterpaper]{amsart}
\usepackage{amsmath,amssymb,amsfonts}
\usepackage{bbm}
\usepackage{enumerate}
\usepackage{amscd, amssymb, latexsym, amsmath, amscd}
\usepackage[all]{xy}
\usepackage{pb-diagram}
\usepackage{mathrsfs}
\usepackage{pictexwd,dcpic}
\usepackage{hyperref,url}
\usepackage{color}

\usepackage{pb-diagram}
\theoremstyle{plain}
\newtheorem{Thm}[equation]{Theorem}

\newtheorem{Cor}[equation]{Corollary}
\newtheorem{Prop}[equation]{Proposition}
\newtheorem{Lem}[equation]{Lemma}
\numberwithin{equation}{section}

\theoremstyle{remark}
\newtheorem{Def}[equation]{Definition}
\newtheorem{Rmk}[equation]{Remark}

\newcommand{\Hom}{\operatorname{Hom}}
\newcommand{\Ext}{\operatorname{Ext}}
\newcommand{\EP}{\operatorname{EP}}
\newcommand{\Gal}{\operatorname{Gal}}
\newcommand{\GL}{\operatorname{GL}}

\newcommand{\SL}{\operatorname{SL}}

\newcommand{\Ind}{\operatorname{Ind}}
\newcommand{\ind}{\operatorname{ind}}
\newcommand{\Jac}{\operatorname{Jac}}

\newcommand{\Irr}{\operatorname{Irr}}

\newcommand{\id}{\operatorname{id}}

\newcommand{\Q}{\mathbb{Q}}

\newcommand{\Z}{\mathbb{Z}}
\newcommand{\R}{\mathbb{R}}

\newcommand{\bm}{\begin{multline*}}
\newcommand{\tu}{\end  {multline*}}

\def\calH{\mathcal{H}}

\def\lra{\longrightarrow}

\title[Ext-GGP model vanishes for tempered representations]
{Ext-analogues of Gan--Gross--Prasad models vanish for tempered representations}
\author{Rui Chen} 

\address{Institute for Advanced Study in Mathematics, Zhejiang University, East No.7 Building, Zijingang Campus, Hangzhou 310058, China }
\email{e0046839@u.nus.edu}
\begin{document}

\begin{abstract}
In this paper we prove that Ext-analogues of Gan--Gross--Prasad models vanish for tempered representations, as conjectured by D. Prasad. In particular, this implies a conjectural Euler--Poincar\'e characteristic formula for Gan--Gross--Prasad models. Besides, we also single out a special class of tempered representations, which includes the Steinberg representation, and show that the restrictions of them are projective. 
\end{abstract}

\maketitle

\textit{Keywords}: Branching laws, Gan--Gross--Prasad conjecture, Ext spaces\\

\textit{2020 Mathematics subject classification}: 11F70, 22E50

\section{Problem and set up}
Let $F$ be a non-Archimedean local field of characteristic $0$ and $E$ either $F$ itself or a quadratic field extension of $F$. Let $V=V_{n+1}$ be an $n+1$-dimensional Hermitian space over $E$, and $W=W_n$ an $n$-dimensional non degenerate subspace of $V$. We denote by $G=G(V)$ and $H=H(W)$ the isometry groups of $V$ and $W$ (so $G$ and $H$ are orthogonal groups or unitary groups), and regard $H$ as a subgroup of $G$ fixing some anisotropic vector. For an irreducible representation $\pi$ of $G$ and $\sigma$ of $H$, one can define the Gan--Gross--Prasad model (``GGP model'' for short)
\[
    \Hom_H(\pi,\sigma).
\]
The study of this model has led to many important results, for example the local Gan--Gross--Prasad conjecture. About ten years ago, Prasad pointed out that, we should not only study the Hom space, but also the Ext-analogue of the GGP model  
\[
    \Ext^i_H(\pi,\sigma),
\]
and the Euler--Poincar\'e characteristic
\[
    \EP_H(\pi,\sigma) = \sum_i (-1)^i\dim \Ext^i_H(\pi,\sigma).
\]
It has been shown by Prasad himself \cite[Thm. 6.1]{MR3966813} that these Ext spaces are all finite dimensional, so $\EP_H(\pi,\sigma)$ is well defined. He also verified that $\EP_H(\pi,\sigma)$ enjoys some good properties, such as the compatibility with parabolic inductions. Comparing with the Riemann--Roch theorem, he conjectured that the Waldspurger's integral formula gives the Euler--Poincar\'e characteristic \cite[Conj. 7.1]{MR3966813}. 

\vskip 5pt

As a first test of Prasad's conjecture, one can consider the GGP model (and its Ext-analogue) for tempered representations. In this case, the Waldspurger's integral formula gives the multiplicity of the GGP model \cite{MR2684300} \cite{MR3202558} \cite{MR4146145}, and Prasad predicted in \cite[Conj. 7.1(2)]{MR3966813} that higher Ext are all vanishing. The main purpose of this paper is to verify this:

\begin{Thm}\label{T:Main}
If both $\pi$ and $\sigma$ are tempered, then 
\[
    \Ext^i_H(\pi,\sigma) = 0
\]
for all $i\geq 1$.
\end{Thm}

\vskip 5pt

%\newpage

\begin{Rmk}
\,
\begin{enumerate}
\item Recently in his IH\'ES lecture note \cite[Thm. 9.3]{prasad2023homological}, Prasad showed that Theorem \ref{T:Main} (see also Remark \ref{R:Highcodim}) implies the same Ext-vanishing result for standard modules (hence also generic representations), as well as the conjectural formula for the Euler--Poincar\'e characteristic \cite[Conj. 7.1]{MR3966813}.

\vskip 5pt

\item Recall that by Casselman's criterion, an irreducible representation is tempered if and only if all exponents of its Jacquet modules are non negative. Therefore the ``temperedness'' can be regarded as a positivity condition. From this point of view, Prasad pointed out that this Ext-vanishing result is an analogue of Kodaira's vanishing theorem in branching laws.

\vskip 5pt

\item One can also consider the GGP model for general linear groups. Similar Ext-vanishing result for general linear groups has been proved by Chan--Savin \cite[Thm. 4.1]{MR4291425} using Bernstein--Zelevinsky derivatives. A key observation in Chan--Savin's work is that one can obtain more information by considering both left and right derivatives at the same time. In our case, the MVW-contragredient trick of Atobe--Gan \cite[Lem. 2.2]{MR3714507} plays the role of this observation.
\end{enumerate}
\end{Rmk}

\vskip 5pt

The main tool that we will use is Atobe's work on Jacquet modules \cite{MR4055180}. In that paper Atobe systematically studied the interaction of Jacquet modules and the local Langlands correspondence (``LLC'' for short) for classical groups, based on some previous results of M{\oe}glin, Jantzen and B. Xu. Using Atobe's result, we can embed a (non cuspidal) tempered representation into a small standard module, and then appeal to the Mackey theory and standard degree shifting argument to study the Hom and Ext spaces. Combining with the aforementioned MVW-contragredient trick, we can eventually show that these Ext spaces are periodic, hence zero. 

\vskip 5pt

Another closely related and interesting question is to determine those (non cuspidal) representations of $G$ whose restrictions to $H$ are projective. In this paper we also investigate a special class of tempered representations (see Definition \ref{D:GSteinberg}). Most of these representations are discrete series, and in some sense, they are tempered representations whose Jacquet modules have maximal possible exponents.
Perhaps the most interesting example among these representations is the Steinberg representation ${\rm St}_G$, defined as 
\[
    {\rm St}_G = C_c^\infty\left(G/P_{min}\right)\big/\sum_{P_{min}\subsetneq Q} C_c^\infty\left(G/Q\right),
\]
where $P_{min}$ is a minimal parabolic subgroup of $G$, and $Q$ runs over all parabolic subgroups of $G$ containing $P_{min}$. Except for the case when $G$ is an even orthogonal group, the L-parameter of ${\rm St}_G$ is irreducible. Indeed, any discrete series representation of $G$ with irreducible L-parameter also belongs to this special class. We refer readers to the list after Definition \ref{D:GSteinberg} for more examples. %functorial lifting of $\pi$ to the general linear group (through the standard representation of the Langlands dual group) is still a discrete series.
Our ingredient for the proof of the projectivity is a uniform lower bound of the largest exponent of Jacquet modules of an arbitrary representation obtained by using the Langlands classification and the Zelevinsky--Aubert duality. Combining this bound with our main result Theorem \ref{T:Main}, we prove the following byproduct:

\begin{Thm}\label{T:ProjSteinberg}
Let $\pi$ be an irreducible tempered representation of $G$, which is of maximal derivative as in Definition \ref{D:GSteinberg}.% and further assume that $\pi$ is a discrete series if $G$ is an odd orthogonal group. 
Then its restriction $\pi~\big|_H$ is projective as a smooth $H$-representation.  
\end{Thm}

\vskip 5pt

As a direct consequence, we obtain:

\begin{Cor}
The following representations of $G$, when restricted to the subgroup $H$, are projective as smooth $H$-representations:

\begin{itemize}
    \item the Steinberg representation ${\rm St}_G$;

    \vskip 5pt

    \item discrete series representation with irreducible L-parameter.
\end{itemize}
\end{Cor}

\vskip 5pt

\begin{Rmk}
\,
\begin{enumerate}
    \item Similar projectivity result for general linear groups has been proved by Chan--Savin \cite[Thm. 5.7]{MR4291425} using Hecke algebra methods. Our approach also works in the general linear group case.  

    \vskip 5pt

\item In \cite{MR4069188}, Chan--Savin also studied the Iwahori component of the Bessel model and obtained some projectivity results for the Steinberg representation of special orthogonal groups.
\end{enumerate}

\end{Rmk}

\vskip 5pt

We end up this section by setting some notations. Let 
\[
    c=\begin{cases}
        \id \quad & \textit{if $E=F$};\\
        \textit{the non-trivial element in }\Gal(E/F) \quad & \textit{if $E\ne F$}. 
    \end{cases}
\]
We denote by
\[
    WD_E=W_E\times \SL_2
\]
the Weil-Deligne group of $E$. An irreducible representation of $WD_E$ will be typically written as $\rho S_a$, where $\rho$ stands for an irreducible representation of $W_E$, and $S_a$ stands for the unique $a$-dimensional irreducible representation of $\SL_2$. Let $\nu$ be the absolute value of $E^\times$, and also regarded as a character of each $\GL_k = \GL_{k/E}$ through the determinant map. We shall use standard notations
\[
    \tau_1\times \tau_2 \quad \text{and} \quad \tau\rtimes\pi
\]
to denote the normalized parabolic induction of general linear groups and classical groups. For a group $G$, its modulus character will be denoted by $\Delta_G$. The hyperbolic plane over $E$ will be denoted by $\calH$.

\vskip 10pt

\section{Preparations for the proofs}

In this section we recall some preliminaries that will be used later. 

\vskip 5pt

\subsection{Jacquet modules}
Let $\rho$ be an irreducible unitary supercuspidal representation of $\GL_k$. By the LLC for general linear groups, $\rho$ corresponds to an irreducible $k$-dimensional representation of $W_E$. We shall often use the same notation to denote $\rho$ and its L-parameter. If $\{x,x-1,\cdots,y\}$ is a segment, we denote by 
\[
    \delta_\rho(x,y) = \langle\rho;x,x-1\cdots,y\rangle
\]
the unique irreducible subrepresentation of $\rho\nu^x\times\rho\nu^{x-1}\times\cdots\times\rho\nu^y$. This is called the generalized Steinberg representation (of the general linear group). 

\vskip 5pt
%In particular, if $x\geq y$, then this is a Steinberg representation and will be also denoted by $\delta_\rho(x,y)$; if $x<y$, then this is a Speh representation and will be also denoted by $\zeta_\rho(x,y)$.\\

Let $\pi$ be a finite length representation of $G=G(V)$. For an irreducible supercuspidal representation $\rho$ of some general linear group, the partial Jacquet module $\Jac_\rho\pi$ is defined as follows. Let $k$ be the dimension of $\rho$ (regarded as an L-parameter), and $P$ the maximal parabolic subgroup of $G$ with Levi component $\GL_k\times G_0$. Here $G_0=G(V_0)$ is the isometry group of a non degenerate subspace $V_0\subset V$ such that $V\simeq V_0 +\calH^k$. If the semi-simplified Jacquet module 
\[
    s.s.\Jac_P\pi = \sum_i\tau_i\boxtimes\sigma_i
\]
for some irreducible representations $\tau_i$ and $\sigma_i$ of $\GL_k$ and $G_0$, then we define the partial Jacquet module as
\[
    \Jac_\rho\pi = \sum_{\{i:\tau_i\simeq \rho\}} \sigma_i.
\]
When $G$ does not have such standard parabolic subgroup $P$, we interpret $\Jac_\rho\pi$ to be $0$. For a string of real numbers
\[
    \underline{y}=(y_1,\cdots,y_r),
\]
we also denote
\[
    \Jac_{\rho,\underline{y}}=\Jac_{\rho\nu^{y_r}}\circ\cdots\circ\Jac_{\rho\nu^{y_1}}.
\]
To simplify notations, when $\rho$ is clear from the context, we omit it from the subscript. The following lemma is worth noting:

\begin{Lem}\label{L:TF}
Let $\rho$ and $\tau$ be two irreducible supercuspidal representations of general linear groups, such that $\rho\not\in\left\{\tau,\left(\tau^c\right)^\vee\right\}$. Then we have
\[
    \Jac_{\rho}(\tau\rtimes\pi) = \tau\rtimes(\Jac_{\rho}\pi)
\]
up to semi-simplifications. %{\color{red}\textit{Check again! Especially for unitary groups there is an issue of taking conjugate, check Tadic's paper...}}
\end{Lem}

\vskip 5pt

\begin{proof}
This is a simple application of the geometric lemma (or, Tadi\'c's formula).

\end{proof}

\vskip 5pt

Since we mainly consider tempered representations, the following version of the Casselman's criterion is also useful to us.

\begin{Lem}\label{L:CasselmanCri}
Let $\pi$ be an irreducible tempered (resp. discrete series) representation of $G(V)$, $V_0\subset V$ a non degenerate subspace such that
\[
    V\simeq V_0+\calH^k,
\]
and $P$ the standard maximal parabolic subgroup of $G$ with Levi component $\GL_k\times G(V_0)$. Then the central exponents of $\Jac_P\pi$, as unramified characters of $\GL_k$, are of the form $\nu^x$ with $x\geq 0$ (resp. $x>0$).

\end{Lem}

\vskip 5pt

\subsection{Local Langlands correspondence}
Now we briefly review the LLC for orthogonal and unitary groups. Readers may consult \cite{MR3135650}, \cite{MR3338302}, \cite{kaletha2014endoscopic}, \cite{MR3708200} and \cite{CZlocal} for more details. Let $G=G(V)$ be the isometry group of the Hermitian space $V$. The LLC for $G$ asserts that there is a finite to one surjective map
\[
    LL: {\Irr}(G) \longrightarrow \Phi(G),
\]
where $\Phi(G)$ is the set of L-parameters of $G$. This map $LL$ preserves various properties such as discreteness and temperedness. 
We recall that an element $M \in \Phi(G)$ is an isomorphism class of conjugate self-dual representation of the Weil-Deligne group $WD_E$ of dimension 
\[
    d_G = \begin{cases}
    \dim V -1 \quad & \text{if $E=F$ and $\dim V$ is odd};\\
    
    \dim V \quad & \text{otherwise}
    \end{cases}
\] 
and sign
\[
    \begin{cases}
    (-1)^{\dim V} \quad & \text{if }E=F;\\
    
    (-1)^{\dim V-1} \quad & \text{if } E\neq F.
    \end{cases}
\] 
When $E=F$ (so $V$ is a quadratic space), there is one more condition that $\det(M)=\chi_V$, the quadratic character of $W_F^{ab}\simeq F^\times$ associated to $V$. Such an L-parameter $M$ is discrete if and only if it is a multiplicity-free direct sum 
\begin{equation*} 
M = \sum_{i} M_i  
\end{equation*}
of irreducible representations $M_i$, such that each $M_i$ is conjugate self-dual of the same sign as $M$; it is tempered if and only if the image of $W_E$ is bounded. For $M\in \Phi(G)$, we put $\Pi_M(G)=LL^{-1}(M)$: this is the local L-packet associated to the L-parameter $M$. 

\vskip 5pt

The following property of the LLC is very important to us. Suppose that $M$ is a tempered L-parameter for $G$, and it has a decomposition
\[
    M = N + M_0 + \left(N^c\right)^\vee,
\]
where $M_0$ and $N$ are representations of $WD_E$. Let $r_V$ be the Witt index of $V$ and $k=\dim N$. We have:
\begin{itemize}
\item if $r_V\geq k$, then 
\[
    \Pi_M(G) = \left\{\pi\subset\tau\rtimes \pi_0~\big|~ \pi_0\in\Pi_{M_0}(G_0)\right\},
\] 
where $G_0=G(V_0)$ for some subspace $V_0\subset V$ such that $V=V_0+\calH^{k}$, and $\tau$ is the irreducible tempered representation of $\GL_k$ with L-parameter $N$;

\vskip 5pt

\item if $r_V< k$, then the L-packet $\Pi_M(G)$ is empty.
\end{itemize}
\vskip 5pt
This property is a part of the local intertwining relation. We remark that the LLC also gives a parametrization inside each L-packet $\Pi_M(G)$, but we will not use this parametrization.

\vskip 5pt

Inspiring by the property \cite[Prop. 9.5]{MR584084} of generalized Steinberg representations of general linear groups, we single out a special class of tempered representations of $G$ as below.  

\begin{Def}\label{D:GSteinberg}
Let $\pi$ be an irreducible tempered representation of $G$. We say that $\pi$ is of maximal derivative, if there exists an irreducible unitary supercuspidal representation $\rho$ of a general linear group $\GL_k$, together with an integer
\[
    a \geq \frac{\dim V}{k}-1,
\]
such that for any irreducible unitary supercuspidal representation $\rho'$ of general linear groups, and any real number $x$, we have $\Jac_{\rho',x}\pi = 0$ unless $\rho'\in \left\{\rho,\left(\rho^c\right)^\vee\right\}$, and $x=\frac{a-1}{2}$.
\end{Def}

\vskip 5pt

As examples, we would like to highlight three families of representations of maximal derivative. We are particularly interested in the last two families.

\begin{enumerate}
    \item Supercuspidal representations: one can simply take arbitrary $\rho$ and arbitrary $a\geq \frac{\dim V}{k}-1$, the requirement in the definition holds automatically.

    \vskip 5pt

    \item Discrete series with irreducible L-parameter $\rho S_a$: for such a discrete series representation $\pi$, the requirement in the definition follows from \cite[Lem. 4.1]{MR4055180}. When $G$ is an odd orthogonal group or a unitary group, the Steinberg representation ${\rm St}_G$ of $G$ has L-parameter $S_{d_G}$ according to \cite[Lem. 1.8]{cunningham2024koszul}, hence belongs to this family.

    \vskip 5pt

    \item Tempered representations with L-parameter of the form $\rho S_a+\delta$, where $\delta$ is irreducible as a representation of $W_E$ and $\dim \delta\leq \dim \rho$ (in the case that $G$ is an odd orthogonal group, we assume that $\dim\delta<\dim \rho$): if $a>1$ or $\delta\not\simeq \left(\rho^c\right)^\vee$, such a representation is a discrete series and is of maximal derivative due to the same reason as in the previous case; if $a=1$ and $\delta\simeq \left(\rho^c\right)^\vee$, one can check by hand. Indeed, if $\pi$ is an irreducible representation of maximal derivative and is not a discrete series, then the assumption $a\geq \frac{\dim V}{k}-1$ forces that $LL\left(\pi\right)=\rho+ \left(\rho^c\right)^\vee$. Again by \cite[Lem. 1.8]{cunningham2024koszul}, when $G$ is an even orthogonal group, the Steinberg representation ${\rm St}_G$ of $G$ has L-parameter $\mathbbm{1}+S_{d_G-1}$, hence belongs to this family.
\end{enumerate} 

\vskip 5pt

\subsection{Role switching}

Now let $\pi$ be an irreducible tempered representation of $G$, and $\sigma$ an irreducible tempered representation of $H$. An important observation is that the role of $(G,\pi)$ and $(H,\sigma)$ are ``symmetric'' in some sense.

\begin{Lem}\label{L:Reduction}
Suppose that
\[
    \pi\simeq \chi\rtimes\pi_0,
\]
where:
\begin{itemize}
\item $V_0\subset V$ is a non degenerate subspace such that $V\simeq V_0+\calH$, and $G_0=G(V_0)$;

\vskip 5pt

\item $P$ the standard maximal parabolic subgroup of $G$ with Levi component $\GL_1\times G_0$;

\vskip 5pt

\item $\chi$ is a unitary character of $\GL_1$, and $\pi_0$ is a tempered representation of $G_0$.
\end{itemize}
\vskip 5pt
Fix an embedding $V_0\hookrightarrow W$. Then we have
\[
    \Ext^i_H\left(\pi,\sigma\right) \simeq \Ext^i_{G_0}\left(\sigma^\vee,\pi_0^\vee\right)
\]
for all $i\geq 0$.
\end{Lem}

\vskip 5pt

\begin{proof}
By the Mackey theory, we know that there is a filtration on $(\chi\rtimes \pi_0)~\big|_H$, whose successive quotients are parametrized by the $H$-orbits of $G/P$. There are at most two such orbits:
\begin{itemize}
    \item the closed orbit consists of isotropic lines of $V$ which are contained in $W$;
    \vskip 5pt

    \item the open orbit consists of isotropic lines of $V$ which are not contained in $W$.
\end{itemize}
\vskip 5pt
Let $\pi_c$ be the quotient of $\chi\rtimes \pi_0$ corresponding to the closed orbit. Then $\pi_c$ can be computed as follows. Let $Q$ be the stabilizer of a representative in the closed orbit. We know that $Q$ is a maximal parabolic subgroup of $H$ stabilizing an isotropic line of $W$. Let $W_0$ be a subspace of $W\cap V_0$ such that $W\simeq W_0+\calH$, and set $H_0=H(W_0)$. Then the Levi component of $Q$ is $M_Q\simeq\GL_1\times H_0$. We have
\begin{align*}
    \pi_c =& \Ind_Q^H\Big(\Delta_P^{\frac{1}{2}}\cdot\chi\boxtimes\left(\pi_0~\big|_{H_0}\right)\Big)\\
    \simeq & \chi\nu^{\frac{1}{2}}\rtimes\left(\pi_0~\big|_{H_0}\right).
\end{align*}
Here we use the symbol ``$\Ind$'' to denote the un-normalized induction. On the other hand, let $\pi_{op}$ be the $H$-subrepresentation of $\chi\rtimes \pi_0$ corresponding to the open orbit. Similar to \cite[Sect. 15]{MR3202556}, we have 
\[
    \pi_{op} \simeq \ind_{G_0}^H \pi_0.
\]
The Mackey theory then yields an $H$-equivariant short exact sequence
\[
    0 \lra \pi_{op} \lra \chi\rtimes \pi_0 \lra \pi_c \lra 0.
\]
Applying the functor $\Hom_H(-, \sigma)$ to this short exact sequence we obtain the long exact sequence
\[
    \cdots \lra \Ext_H^i\left(\pi_c, \sigma\right) \lra \Ext_H^i\left(\pi,\sigma\right) \lra \Ext_H^i\left(\pi_{op},\sigma\right) \lra \Ext_H^{i+1}\left(\pi_c,\sigma\right) \lra \cdots
\]
Since $\sigma$ is tempered, by Lemma \ref{L:CasselmanCri} any irreducible constituent of the normalized Jacquet module $\Jac_{\overline{Q}}\sigma$ has non positive exponents. Hence
\[
    \Ext_H^i\left(\pi_c, \sigma\right) = \Ext_{M_Q}^i\left(\chi\nu^{\frac{1}{2}}\boxtimes\left(\pi_0~\big|_{H_0}\right) , \Jac_{\overline{Q}}\sigma\right)
\]
vanishes for all degree $i$. The long exact sequence then gives 
\begin{align*}
    \Ext_H^i\left(\pi,\sigma\right) &\simeq  \Ext_H^i\left(\pi_{op},\sigma\right)\simeq \Ext_{G_0}^i\left(\sigma^\vee, \pi_0^\vee\right).
\end{align*}
This completes the proof.

\end{proof}

\vskip 5pt

As a corollary, fix an embedding $V\hookrightarrow W+\calH$ and apply this lemma, we get:
\begin{Cor}\label{C:Switch}
For any non conjugate self-dual unitary character $\chi$ of $\GL_1$ we have 
\begin{equation*}
    \Ext_H^i(\pi,\sigma) \simeq \Ext_G^i\left(\chi\rtimes \sigma^\vee, \pi^\vee\right)
\end{equation*}
for all $i\geq 0$.
\end{Cor}

\vskip 10pt

\section{Proof of the main theorem I: discrete series case}

Now we are ready to prove our main result: Theorem \ref{T:Main}. In this section we aim at the discrete series case, but for the convenience of using induction argument, we still allow the L-parameters containing certain bad parity part.

\subsection{Working assumption}\label{SS:WA}
Let $M$ be a tempered L-parameter of $G$ and $N$ of $H$. Write $M$ as 
\[
    M= M_{np} + M_{bp} + (M_{np}^c)^\vee,
\]
where: 
\begin{itemize}
\item $M_{bp}$ is the summation of all irreducible constituents which are conjugate self-dual of the same parity as $M$;
\vskip 5pt

\item $M_{np}$ is the bad parity part of $M$.
\end{itemize}
Assume that:

\begin{itemize}
    \item[(\textbf{A.0})] $M_{bp}$ is a discrete L-parameter, and 
    \[
    M_{np} = m_1\chi_1 +\cdots + m_r\chi_r,
\]
where for each $1\leq i\leq r$, $\chi_i$ is a non conjugate self-dual character of $W_E$ with multiplicity $m_i$, such that for any $i\neq j$, we have
\[
    \chi_i \neq \chi_j \quad \text{and} \quad \chi_i \neq \left(\chi_j^c\right)^\vee.
\]
    %$M_{np}$ is a summation of non conjugate self-dual characters of $W_E$.
\end{itemize}
\vskip 5pt
Note that any discrete L-parameter satisfies this assumption. We also impose the same condition on $N$. In this section we show that:
\begin{Prop}\label{P:discCase}
Suppose that $M$ and $N$ satisfy (\textbf{A.0}). Then for any $\pi\in\Pi_M(G)$ and $\sigma\in\Pi_N(H)$, we have
\[
    \Ext^i_H(\pi,\sigma) = 0
\]
for all $i\geq 1$.
\end{Prop}

\vskip 5pt

We shall prove this proposition by induction on $\dim M_{bp} + \dim N_{bp}$. During the proof we will use the same argument as \cite[Thm. 15.1]{MR3202556} many times. The following lemma guarantees the legitimacy of our arguments.

\begin{Lem}
There are infinitely many (non conjugate self-dual) unitary characters of $\GL_1$ up to unramified twists.
\end{Lem}

\vskip 5pt

\begin{proof}
Easy to see from the isomorphism
\[
E^\times \simeq\varpi^\Z\times k_E^\times \times \mu_{p^\infty}(E) \times \Z_p^d.
\]
Here $\varpi$ is a uniformizer, $k_E$ is the residue field of $E$, $\mu_{p^\infty}(E)$ is the group of roots of unity of $p$-power order in $E$, and $d$ is the degree of $E$ over $\Q_p$. 

\end{proof}

\vskip 5pt

\subsection{Some small standard modules}

Given an L-parameter $M$ satisfies (\textbf{A.0}). By the property of the LLC, we know that any $\pi\in\Pi_M(G)$ can be written as an irreducible parabolic induction
\[
    \pi \simeq \chi_1^{m_1}\times \cdots \times\chi_r^{m_r} \rtimes \pi_{bp},
\]
where $\pi_{bp}\in\Pi_{M_{bp}}(G(V_{bp}))$, with $V_{bp}\subset V$ a non degenerate subspace such that $V\simeq V_{bp}+\calH^{m}$, and $m=m_1+\cdots+m_r$. 

\vskip 5pt

Suppose that there is some irreducible unitary supercuspidal representation $\rho$ of a general linear group, such that
\[
    \Jac_{\rho,x}\pi_{bp} \neq 0
\]
for some $x\in\R$. Then by \cite[Lem. 4.1]{MR4055180} we have $x\in\frac{1}{2}\Z_{>0}$, and $\rho S_{2x+1}\subset M_{bp}$. The following two lemmas, which were originally stated in Atobe's paper for discrete (respectively, good parity) L-parameters, are key ingredients of our proof of the main theorem.

\begin{Lem}\label{L:StdMod-1}
In the context of above discussion, suppose that $x=\frac{1}{2}$, or $\rho S_{2x-1}\not\subset M_{bp}$. Let $V_0\subset V$ be a non degenerate subspace such that $V\simeq V_0 +\calH^k$, where $k$ is the dimension of $\rho$ regarded as an L-parameter. Set $G_0=G(V_0)$. Then there is a short exact sequence %{\color{red}\textit{Atobe's result are for discrete/ good partity L-parameters! Need to generalize slightly...}}
\begin{equation*}
    0 \lra \pi \lra \rho\nu^x\rtimes \pi_0 \lra \pi' \lra 0,
\end{equation*}
where $\pi_0\in\Pi_{M_0}(G_0)$ with 
\[
    M_0 = M-\rho S_{2x+1} + \rho S_{2x-1},
\]
and $\pi'$ is the unique irreducible quotient of $\rho\nu^x\rtimes \pi_0$. 
\end{Lem}

\vskip 5pt

\begin{proof}
We prove this lemma by induction on $m=\dim M_{np}$. The basic case is that $M_{np}=0$, i.e. $M= M_{bp}$ is a discrete L-parameter. Then the assertion has been proved in \cite[Lem. 5.1]{MR4055180}. 

\vskip 5pt

Now suppose that we have proved the lemma for any L-parameter $M'$ satisfying our requirements and $\dim M'_{np}<m$. Let
\[
    M_1 = M - m_1\chi_1 - m_1\left(\chi_1^c\right)^\vee,
\]
and also
\[
    \pi_1 = \chi_2^{m_2}\times\cdots\times\chi_r^{m_r}\rtimes \pi_{bp}.
\]
Then the L-parameter associated to $\pi_1$ is $M_1$, and $\pi=\chi_1^{m_1}\rtimes \pi_1$. By the induction hypothesis, we have a short exact sequence
\begin{equation*}%\label{E:S.E.S-0}\tag{$\natural$}
    0 \lra {\pi_1} \lra \rho\nu^x\rtimes {\pi_{1,0}} \lra {\pi'_1} \lra 0,
\end{equation*}
where $\pi_{1,0}$ is an irreducible tempered representation with L-parameter
\[
    M_{1,0} = M_1-\rho S_{2x+1} + \rho S_{2x-1},
\]
and $\pi'_1$ is the unique irreducible quotient of $\rho\nu^x\rtimes {\pi_{1,0}}$. Applying the functor $\chi_1^{m_1}\rtimes-$ to this short exact sequence, we get
\[
    0 \lra \pi \lra \chi_1^{m_1}\times\rho\nu^x\rtimes {\pi_{1,0}} \lra \chi_1^{m_1}\rtimes{\pi'_1} \lra 0.
\]
Since $\chi_1\neq \rho$, we have
\[
    \chi_1^{m_1}\times\rho\nu^x\rtimes {\pi_{1,0}} \simeq \rho\nu^x\rtimes\left(\chi_1^{m_1}\rtimes {\pi_{1,0}}\right).
\]
By the property of the LLC, $\chi_1^{m_1}\rtimes {\pi_{1,0}}$ is irreducible tempered with L-parameter $M_0$. Therefore to complete the proof it only remains to show that $\chi_1^{m_1}\rtimes{\pi'_1}$ is irreducible.

\vskip 5pt

To show that $\chi_1^{m_1}\rtimes{\pi'_1}$ is irreducible, firstly we claim that $\pi_1'$ is $\chi_1$-reduced, namely
\[
    \Jac_{\chi_1}\pi_1' = 0.
\]
Indeed, by Lemma \ref{L:TF} and the property of the LLC we have
\[
    \Jac_{\chi_1}\left(\rho\nu^x\rtimes {\pi_{1,0}}\right) = \rho\nu^x\times\chi_2^{m_2}\times\cdots\times\chi_r^{m_r}\rtimes \Jac_{\chi_1}{\pi_{0,bp}},
\]
where $\pi_{0,bp}$ is an irreducible tempered representation with L-parameter 
\[
    M_{bp} -\rho S_{2x+1} + \rho S_{2x-1}.
\]
Since $\chi_1\not\subset M_{bp}$, it follows from \cite[Lem. 4.1]{MR4055180} that $\Jac_{\chi_1}{\pi_{0,bp}}=0$. This implies our claim. Let $\pi'$ be the unique irreducible subrepresentation of $\chi_1^{m_1}\rtimes{\pi'_1}$. We have
\[
    \left(\Jac_{\chi_1}\right)^{m_1} \pi' \neq 0,
\]
and $\pi'$ is the only subquotient of $\chi_1^{m_1}\rtimes{\pi'_1}$ with this property. On the other hand, by the property of the LLC, we have
\[
    \chi_1^{m_1}\rtimes{\pi_{1,0}} \simeq \left(\left(\chi_1^c\right)^\vee\right)^{m_1}\rtimes{\pi_{1,0}}.
\]
Let $\pi''$ be the unique irreducible quotient of $\chi_1^{m_1}\rtimes{\pi'_1}$. Then $\pi''$ is also a quotient of 
\[
    \rho\nu^x\rtimes\left(\chi_1^{m_1}\rtimes {\pi_{1,0}}\right) \simeq \left(\left(\chi_1^c\right)^\vee\right)^{m_1} \rtimes \left(\rho\nu^x\rtimes {\pi_{1,0}}\right),
\] 
which implies that
\[
    \left(\Jac_{\chi_1}\right)^{m_1} \pi'' \neq 0.
\]
Thus $\pi' = \pi''$, i.e. $\chi_1^{m_1}\rtimes{\pi'_1}$ is irreducible. This completes the proof.

\end{proof}

\vskip 5pt

\begin{Lem}\label{L:StdMod-2}
In the context of above discussion, suppose that $x>\frac{1}{2}$ and $\rho S_{2x-1}\subset M_{bp}$. Let $V_0\subset V$ be a non degenerate subspace such that $V\simeq V_0 +\calH^{k\cdot 2x}$, where $k$ is the dimension of $\rho$ regarded as an L-parameter. Set $G_0=G(V_0)$. Then there is a short exact sequence %{\color{red}\textit{Atobe's result are for discrete/ good partity L-parameters! Need to generalize slightly...}}
\begin{equation*}%\label{E:S.E.S-3}\tag{$\clubsuit$}
    0 \lra \pi_1\oplus \pi_2 \lra \delta_\rho(x, -x+1)\rtimes \pi_0 \lra \pi' \lra 0,
\end{equation*}
where $\pi_0\in\Pi_{M_0}(G_0)$ with L-parameter
\[
    M_0 = M-\rho S_{2x+1} - \rho S_{2x-1},
\]
both $\pi_1$ and $\pi_2$ belong to $\Pi_M(G)$, $\pi\in\{\pi_1,\pi_2\}$ and $\pi'$ is the unique irreducible quotient of $\delta_\rho(x, -x+1)\rtimes \pi_0$. 
\end{Lem}

\vskip 5pt

\begin{proof}
Similarly, in the case that $M=M_{bp}$ is of good parity, the assertion has been proved in \cite[Prop. 5.2]{MR4055180}; the general case can be shown by induction on $\dim M_{np}$.

\end{proof}

\vskip 5pt

\subsection{Basic case}

Back to the proof of Proposition \ref{P:discCase}, we first show the basic case, namely the case that
\[
    \dim M_{bp} + \dim N_{bp}= \dim V_{an} + \dim W_{an}.
\]
where $V_{an}$ and $W_{an}$ are the anisotropic kernels of $V$ and $W$. Let 
\[
    M = \chi+M_0+\left(\chi^c\right)^\vee
\]
for some non conjugate self-dual unitary character $\chi$. Then for any $\pi\in \Pi_M(G)$, there exists some $\pi_0\in\Pi_{M_0}(G_0)$ such that
\[
    \pi \simeq \chi\rtimes \pi_0.
\]
Here $G_0$ is the isometry group of a non degenerate space $V_0\subset V$ such that $V\simeq V_0+\calH$. By Lemma \ref{L:Reduction}, we know that
\[
    \Ext^i_H\left(\pi,\sigma\right) \simeq  \Ext^i_{G_0}\left(\sigma^\vee,\pi_0^\vee\right)
\]
for all $i\geq 0$. Doing this repeatedly, finally we can ``peel off'' all of the bad parity part, and our target $\Ext_H^i\left(\pi,\sigma\right)$ is isomorphic to some Ext space over $G(V_{an})$ or $H(W_{an})$. As both $G(V_{an})$ and $H(W_{an})$ are compact, we know that
\[
    \Ext_H^i\left(\pi,\sigma\right) = 0
\]
for all $i\geq 1$. This completes the proof of the basic case.

\vskip 5pt

\subsection{Induction step}\label{SS:Induction}

By the property of the LLC, we know that any $\pi\in\Pi_M(G)$ can be written as an irreducible parabolic induction
\[
    \pi = \tau \rtimes \pi_{bp},
\]
where $\tau$ is an irreducible tempered representation of some general linear group, and $\pi_{bp}\in\Pi_{M_{bp}}(G(V_{bp}))$, where $V_{bp}\subset V$ is a non degenerate subspace of proper dimension. Similarly, we can also write $\sigma\in\Pi_N(H)$ as
\[
    \sigma = \lambda \rtimes \sigma_{bp}
\]
for an irreducible tempered representation $\lambda$ of a general linear group, and $\sigma_{bp}\in\Pi_{N_{bp}}(H(W_{bp}))$.

\vskip 5pt

If both $\pi_{bp}$ and $\sigma_{bp}$ are supercuspidal, then similar to the basic case, by applying Lemma \ref{L:Reduction} one can show that
\[
    \Ext_H^i(\pi,\sigma) \simeq \Ext_{H'}^i\left(\pi',\sigma'\right)
\]
with $\pi'\simeq \pi_{bp}$ or $\sigma'\simeq\sigma_{bp}$ supercuspidal. Hence in this case $\Ext_H^i(\pi,\sigma)$ vanishes for all $i\geq 1$.

\vskip 5pt

So we can assume that at least one of $\pi_{bp}$ and $\sigma_{bp}$ is not supercuspidal. Without lose of the generality we may further assume that:
\vskip 5pt
\begin{itemize}
\item[(\textbf{A.1})] $\pi_{bp}$ is not supercuspidal. 
\end{itemize}
\vskip 5pt
This is because that if $\pi_{bp}$ is supercuspidal but $\sigma_{bp}$ not, by Corollary \ref{C:Switch} we can replace $(\pi,\sigma)$ by $\left(\chi\rtimes \sigma^\vee, \pi^\vee\right)$. 

\vskip 5pt

Given an irreducible unitary supercuspidal representation $\rho'$ of a general linear group, if
\[
    \Jac_{\rho', x'}\pi_{bp} \neq 0
\]
for some $x'\in\R$, then by \cite[Lem. 4.1]{MR4055180} we have $x'\in \frac{1}{2}\Z_{> 0}$, and $\rho' S_{2x'+1}\subset M_{bp}$. Let $x\in\frac{1}{2}\Z_{> 0}$ be the minimal half-integer such that there exists an irreducible unitary supercuspidal representation $\rho$ of a general linear group, with the property that
\[
    \Jac_{\rho, x}\pi_{bp} \neq 0.
\]
Once again we can further assume that:
\vskip 5pt
\begin{itemize}
\item[(\textbf{A.2})] For any irreducible unitary supercuspidal representation $\rho'$ of a general linear group, we have
\[
    \Jac_{\rho', x-\frac{1}{2}}\sigma_{bp} =0.
\]
\end{itemize}
\vskip 5pt
Because otherwise we can replace $(\pi,\sigma)$ by $\left(\chi\rtimes \sigma^\vee, \pi^\vee\right)$ and replace $x$ by $x-\frac{1}{2}$, thanks to Corollary \ref{C:Switch}. We put $a=2x+1 \geq 2$. Finally, combining Corollary \ref{C:Switch} with the LLC, we can further assume that:
\vskip 5pt
\begin{itemize}
\item[(\textbf{A.3})] Any character $\chi$ occurring in $N_{np}$ or $\left(N^c_{np}\right)^\vee$ is not an unramified twist of $\rho$.
\end{itemize}
\vskip 5pt

There are two subcases.

\vskip 10pt
%%%%%%%%%%%%%%%%%%%%%%%%%%%%%%%%%%%%%%%%%%%%%%%%%%%%%%%%%%%%%%%%%%%%%%%%%%%%%%%%%%%
\iffalse$\bullet$ \textit{Subcase 1}: Suppose that $a=1$. Then by \cite[Thm. 4.2(4)]{MR4055180}, we know that $2\rho\subset M_{bp}$. Let $V_0\subset V$ be a non degenerate subspace such that $V\simeq V_0 +\calH^k$, where $k$ is the dimension of $\rho$ regarded as an L-parameter. Set $G_0=G(V_0)$, and $P$ the maximal parabolic subgroup stabilizing a $k$-dimensional isotropic subspace of $V$. It follows from the property of the LLC that $\pi$ is a direct summand of $\rho\rtimes \pi_0$, where $\pi_0\in\Pi_{M_0}(G_0)$ with 
\[
    M_0 = M-2\rho.
\]
Using the same argument as in the basic case (see also \cite[Sect. 15]{MR3202556}), we deduce that 
\[
    \Ext_H^i(\rho\rtimes \pi_0,\sigma) \simeq \Ext_H^i(\pi_{op}, \sigma),
\]
where $\pi_{op}$ is the $H$-subrepresentation of $\rho\rtimes \pi_0$ corresponding to the open $H$-orbit of $G/P$. 

\vskip 10pt
\fi%%%%%%%%%%%%%%%%%%%%%%%%%%%%%%%%%%%%%%%%%%%%%%%%%%%%%%%%%%%%%%%%%%%%%%%%%%%%%%%%%%%%%%%%%%%%%%%%%%%%%%%%%%%%%%%%%%%%

$\bullet$ \textit{Subcase 1}: Suppose that $a=2$, or $\rho S_{a-2}\not\subset M_{bp}$. Let $V_0\subset V$ be a non degenerate subspace such that $V\simeq V_0 +\calH^k$, where $k$ is the dimension of $\rho$ regarded as an L-parameter. Set $G_0=G(V_0)$, and $P$ the maximal parabolic subgroup stabilizing a $k$-dimensional isotropic subspace of $V$. By Lemma \ref{L:StdMod-1}, there is a short exact sequence
%\cite[Lem. 5.1]{MR4055180}, there is a short exact sequence {\color{red}\textit{Atobe's result are for discrete/ good partity L-parameters! Need to generalize slightly...}}
\begin{equation*}\label{E:S.E.S-1}\tag{$\spadesuit$}
    0 \lra \pi \lra \rho\nu^x\rtimes \pi_0 \lra \pi' \lra 0,
\end{equation*}
where $\pi_0\in\Pi_{M_0}(G_0)$ with 
\[
    M_0 = M-\rho S_a + \rho S_{a-2},
\]
and $\pi'$ is the unique irreducible quotient of $\rho\nu^x\rtimes \pi_0$. Firstly we would like to show that
\begin{equation*}\label{E:AuxVanish1}\tag{$\dagger.1$}
    \Ext_H^i\left(\rho\nu^x\rtimes \pi_0, \sigma\right) = 0
\end{equation*}
for all $i\geq 1$. Similar to the basic case, we use the Mackey theory to compute this Ext space. The quotient $\Pi_c$ of $\rho\nu^x\rtimes \pi_0$ corresponding to the closed $H$-orbit of $G/P$ is
\[
    \Pi_c = \rho\nu^{x+\frac{1}{2}}\rtimes \left(\pi_0~\big|_{H_0} \right),
\]
which does not make any contribution to the Ext space since $\sigma$ is tempered and $x+\frac{1}{2}>0$. Similar to \cite[Thm. 15.1]{MR3202556}, the $H$-subrepresentation $\Pi_{op}$ corresponding to the open $H$-orbit of $G/P$ is
\[
    \Pi_{op} = \ind_Q^H\left(\ind_{U_k}^R\psi_k\right)\boxtimes \pi_0,
\]
where $Q$ is the stabilizer of a representative in the open orbit, $R$ is a mirabolic subgroup of $\GL_k$, $U_k$ is the unipotent radical of a Borel subgroup of $\GL_k$, and $\psi_k$ is a generic character of $U_k$. Here we have made use of a well-known result of Gelfand--Kazhdan which asserts that 
\[
   \Delta_P^{\frac{1}{2}}\cdot \rho\nu^x~\big|_R \simeq \ind_{U_k}^R\psi_k.
\]
So we have
\[
    \Ext_H^i\left(\rho\nu^x\rtimes \pi_0, \sigma\right) \simeq \Ext_H^i\left(\Pi_{op}, \sigma\right).
\]
Let $\chi_1,\cdots,\chi_k$ be a sequence of unitary non conjugate self-dual characters of $\GL_1$, such that $\chi_j$ does not occur in the cuspidal component of $\sigma$ for all $j\in\{1,\cdots,k\}$. Set
\[
    \xi = \chi_1 \times \cdots \times \chi_k.
\]
We consider the auxiliary representation $\xi\rtimes \pi_0$. It is an irreducible tempered representation of $G$, with L-parameter
\[
    M^\flat = \left(\chi_1+\cdots+\chi_k\right) + M_0 +\left(\chi_1^c+\cdots+\chi_k^c\right)^\vee.
\]
On the one hand, since $\dim M^\flat_{bp} < \dim M_{bp}$, it follows from the induction hypothesis that
\[
    \Ext_H^i(\xi\rtimes \pi_0,\sigma) = 0.
\]
On the other hand these Ext spaces can be computed using the Mackey theory as well. The quotient $\Pi^\flat_c$ of $\xi\rtimes \pi_0$ corresponding to the closed $H$-orbit of $G/P$ is
\[
    \Pi^\flat_c = \xi\nu^{\frac{1}{2}}\rtimes \left(\pi_0~\big|_{H_0} \right),
\]
which will not contribute to the Ext space since $\sigma$ is tempered. The $H$-subrepresentation $\Pi^\flat_{op}$ corresponding to the open $H$-orbit of $G/P$ is
\[
    \Pi^\flat_{op} = \ind_Q^H\left(\Delta_P^{\frac{1}{2}}\cdot\xi~\big|_R\right)\boxtimes \pi_0,
\]
By a classical result of \cite[Prop. 3.2, 3.5]{MR579172}, the restriction of $\xi$ to $R$ has a filtration with composition factors
\[
    \ind_{T_\alpha\rtimes N_\alpha}^{R}\left(\nu^{\frac{\alpha}{2}}\xi^{(\alpha)}\boxtimes\psi_\alpha\right),
\]
where:
\begin{itemize}
    \item $T_\alpha$ is a subgroup of $R$ isomorphic to $\GL_{k-\alpha}\times U_\alpha$, with $U_\alpha$ the subgroup of strictly upper triangular matrices in $\GL_\alpha$;
    \vskip 5pt

    \item $N_\alpha$ is the unipotent radical of $P_\alpha$, where $P_\alpha$ is the standard parabolic subgroup of $\GL_{k}$ with Levi component $\GL_{k-\alpha}\times \GL_\alpha$;
    \vskip 5pt

    \item $\xi^{(\alpha)}$ is the $\alpha$-th Bernstein--Zelevinsky derivative of $\xi$;
    \vskip 5pt

    \item $\psi_\alpha$ is a generic character of $U_\alpha$.
\end{itemize}
\vskip 5pt
We have
\[
    \Ext^i_H\left(\ind_Q^H\left(\Delta_P^{\frac{1}{2}}\cdot\ind_{T_\alpha\rtimes N_\alpha}^{R}\left(\nu^{\frac{\alpha}{2}}\xi^{(\alpha)}\boxtimes\psi_\alpha\right)\right)\boxtimes \pi_0,\sigma\right) = \Ext^i_H\left(\nu^{\frac{1}{2}}\xi^{(\alpha)}\rtimes\left(\ind_{B_\alpha}^{H_\alpha}\mu_\alpha^\vee\otimes\pi_0\right),\sigma\right),
\]
where $H_\alpha=H(W_\alpha)$ for some subspace $W_\alpha\subset W$ such that $W=W_\alpha+\calH^{k-\alpha}$, $B_\alpha$ is the Bessel subgroup of $H_\alpha$ associated to the pair $(W_\alpha,V_0)$, and $\mu_\alpha$ is the distinguished character attached to $B_\alpha$ (see \cite[Sect. 12]{MR3202556}). 
When $\alpha<k$, $\xi^{(\alpha)}$ is a product of $\chi_j$'s; since we have assumed that $\chi_j$ does not occur in the cuspidal component of $\sigma$, it will not contribute to the Ext space. Thus
\[
    \Ext_H^i\left(\Pi^\flat_{op},\sigma\right) \simeq \Ext_H^i\left(\ind_Q^H\left(\ind_{U_k}^R\psi_k\right)\boxtimes \pi_0,\sigma\right).
\]  
This implies that
\[
    \Ext_H^i\left(\Pi_{op},\sigma\right) \simeq \Ext_H^i\left(\Pi^\flat_{op},\sigma\right)
\]
vanishes for all $i\geq 1$, and hence (\ref{E:AuxVanish1}) follows. Applying the functor $\Hom_H(-,\sigma)$ to the short exact sequence (\ref{E:S.E.S-1}), we deduce that
\[
    \Ext^i_H\left(\pi,\sigma\right) \simeq \Ext^{i+1}_H\left(\pi',\sigma\right)
\]
for all $i \geq 1$.

\vskip 5pt

Dually, applying both the contragredient and MVW-involution to (\ref{E:S.E.S-1}), we get a dualized short exact sequence
\begin{equation*}\label{E:S.E.S-2}\tag{$\clubsuit$}
    0 \lra \pi' \lra \rho\nu^{-x}\rtimes \pi_0 \lra \pi \lra 0,
\end{equation*}
We would also like to show that
\begin{equation*}\label{E:AuxVanish2}\tag{$\dagger.2$}
    \Ext_H^i\left(\rho\nu^{-x}\rtimes \pi_0, \sigma\right) = 0
\end{equation*}
for all $i\geq 1$. Similarly, we appeal to the Mackey theory again. The quotient $\Pi'_c$ of $\rho\nu^{-x}\rtimes \pi_0$ corresponding to the closed orbit is
\[
    \Pi'_c = \rho\nu^{-x+\frac{1}{2}}\rtimes \left(\pi_0~\big|_{H_0} \right).
\]
To see $\Pi'_c$ does not contribute to the Ext spaces, we use the second Frobenius reciprocity
\begin{equation*}%\label{E:Subcase1-Aux1}\tag{$\dagger.3$}
    \Ext^i_{H}\left(\Pi'_c, \sigma\right) = \Ext^i_{\GL_k\times H_0}\left(\rho\nu^{-x+\frac{1}{2}}\boxtimes \left(\pi_0~\big|_{H_0} \right), \Jac_{\overline{Q_c}}\sigma\right),
\end{equation*}
where $Q_c$ is the stabilizer of a representative in the closed orbit, which is also a maximal parabolic subgroup of $H$ with Levi component $\GL_k\times H_0$. Recall that $\sigma = \lambda\rtimes \sigma_{bp}$, the semi-simplification of $\Jac_{\overline{Q_c}}\sigma$ can be analyzed using the Tadi\'c's formula: if we write
\[
    s.s.\Jac_{\overline{Q_c}}\sigma = \sum_{j} \tau_j\boxtimes \sigma_j
\]
for some irreducible representations $\tau_j$ and $\sigma_j$ of $\GL_k$ and $H_0$, then for each $\tau_j$ it has three possibilities:
\begin{itemize}
\item $\tau_j$ is not supercuspidal;

\vskip 5pt

\item $\tau_j$ is supercuspidal but not an unramified twist of $\rho$; 

\vskip 5pt

\item $\tau_j$ is an unramified twist of $\rho$; by the assumption (\textbf{A.3}), we know that $\tau_j$ arises from taking Jacquet modules of $\sigma_{bp}$; then by the assumption (\textbf{A.2}), the central exponent of $\tau_j$ is not equal to $\nu^{-x+\frac{1}{2}}$.
\end{itemize}
\vskip 5pt
In any of these possibilities, we always have
\[
    \Ext^i_{\GL_k}\left(\rho\nu^{-x+\frac{1}{2}} , \tau_j\right) = 0
\]
for all $i\geq 0$. Hence we deduce from the K\"unneth's formula that
\[
    \Ext^i_{H}\left(\Pi'_c, \sigma\right) = 0
\]
for all $i\geq 0$. As for the $H$-subrepresentation $\Pi'_{op}$ corresponding to the open orbit, we have
\[
    \Pi'_{op} \simeq \Pi_{op}.
\]
Hence (\ref{E:AuxVanish2}) follows. Applying the functor $\Hom_H(-,\sigma)$ to the short exact sequence (\ref{E:S.E.S-2}), we deduce that
\[
    \Ext^i_H\left(\pi',\sigma\right) \simeq \Ext^{i+1}_H\left(\pi,\sigma\right)
\]
for all $i \geq 1$. From these isomorphisms one can see that $\Ext^i_H(\pi,\sigma)$ is periodic, i.e.
\[
    \Ext^i_H(\pi,\sigma) \simeq \Ext^{i+2}_H(\pi,\sigma)
\]
for all $i \geq 1$. Since the higher extensions vanish when the degree is sufficiently large \cite[Pg. 98, Sect. 4.2]{BNote}, these groups $\Ext^i_H(\pi,\sigma)$ must vanish for all $i \geq 1$ with no other choice.

\vskip 10pt

$\bullet$ \textit{Subcase 2}: Suppose that $a>2$ and $\rho S_{a-2}\subset M_{bp}$. Let $V_0\subset V$ be a non degenerate subspace such that $V\simeq V_0 +\calH^{k(a-1)}$, where $k$ is the dimension of $\rho$ regarded as an L-parameter. Set $G_0=G(V_0)$, and $P$ the maximal parabolic subgroup stabilizing a $k(a-1)$-dimensional isotropic subspace of $V$. By Lemma \ref{L:StdMod-2}, there is a short exact sequence
%\cite[Prop. 5.2]{MR4055180}, there is a short exact sequence {\color{red}\textit{Atobe's result are for discrete/ good partity L-parameters! Need to generalize slightly...}}
\begin{equation*}\label{E:S.E.S-3}\tag{$\heartsuit$}
    0 \lra \pi_1\oplus \pi_2 \lra \delta_\rho(x, -x+1)\rtimes \pi_0 \lra \pi' \lra 0,
\end{equation*}
where $\pi_0\in\Pi_{M_0}(G_0)$ with L-parameter
\[
    M_0 = M-\rho S_a - \rho S_{a-2},
\]
both $\pi_1$ and $\pi_2$ belong to $\Pi_M(G)$, $\pi\in\{\pi_1,\pi_2\}$ and $\pi'$ is the unique irreducible quotient of $\delta_\rho(x, -x+1)\rtimes \pi_0$. Similar to the previous subcase, we would like to show that
\begin{equation*}\label{E:AuxVanish3}\tag{$\ddagger.1$}
    \Ext_H^i\left(\delta_\rho(x, -x+1)\rtimes \pi_0, \sigma\right) = 0
\end{equation*}
for all $i\geq 1$ using the Mackey theory. The quotient $\Pi_c$ of $\delta_\rho(x, -x+1)\rtimes \pi_0$ corresponding to the closed $H$-orbit of $G/P$ is
\[
    \Pi_c = \delta_\rho\left(x+\frac{1}{2}, -x+\frac{3}{2}\right)\rtimes \left(\pi_0~\big|_{H_0} \right),
\]
which does not make any contribution to the Ext space by the temperedness of $\sigma$. The $H$-subrepresentation $\Pi_{op}$ corresponding to the open $H$-orbit of $G/P$ is 
\[
    \Pi_{op} = \ind_Q^H\left(\Delta_P^{\frac{1}{2}}\cdot\delta_\rho(x, -x+1)~\big|_R\right)\boxtimes \pi_0,
\]
where $Q$ is the stabilizer of a representative in the open orbit, $R$ is a mirabolic subgroup of $\GL_{k(a-1)}$. By \cite[Prop. 3.2, 3.5]{MR579172}, the restriction $\delta_\rho(x, -x+1)~\big|_R$ is ``glued'' from composition factors
\[
    \ind_{T_\alpha\rtimes N_\alpha}^{R}\left(\nu^{\frac{\alpha}{2}}\delta_\rho(x, -x+1)^{(\alpha)}\boxtimes\psi_\alpha\right).
\]
Here $T_\alpha$, $N_\alpha$, $\delta_\rho(x, -x+1)^{(\alpha)}$ and $\psi_\alpha$ are defined similar to Subcase 1. The Bernstein--Zelevinsky derivative $\delta_\rho(x, -x+1)^{(\alpha)}$ is non-zero if and only if $\alpha = k\cdot j$ for some $j\in\{1,\cdots, a-1\}$, and
\[
    \delta_\rho(x, -x+1)^{(k\cdot j)} = \delta_\rho(x, -x+j+1).
\]
When $\alpha = k\cdot j<k(a-1)$, we have
\[
    \Ext^i_H\left(\ind_Q^H\left(\Delta_P^{\frac{1}{2}}\cdot\ind_{T_\alpha\rtimes N_\alpha}^{R}\left(\nu^{\frac{\alpha}{2}}\delta_\rho(x, -x+j+1)\boxtimes\psi_\alpha\right)\right)\boxtimes \pi_0,\sigma\right)
\]
\[
\simeq \Ext^i_H\left(\delta_\rho\left(x+\frac{1}{2}, -x+j+\frac{3}{2}\right)\rtimes\left(\ind_{B_\alpha}^{H_\alpha}\mu_\alpha^\vee\otimes\pi_0\right),\sigma\right)
\]
vanishes by the temperedness of $\sigma$. So
\[
    \Ext_H^i\left(\Pi_{op}, \sigma\right) \simeq \Ext_H^i\left(\ind_Q^H\left(\ind_{U_{k(a-1)}}^R\psi_{k(a-1)}\right)\boxtimes \pi_0,\sigma\right).
\]
Let $\chi_1,\cdots,\chi_{k(a-1)}$ be a sequence of non conjugate self-dual characters of $\GL_1$, and set $\xi=\chi_1\times\cdots\times\chi_{k(a-1)}$. Consider the auxiliary representation $\xi\rtimes\pi_0$. As we have seen before, these characters can be chosen suitably so that
\[
    \Ext^i_H\left(\xi\rtimes\pi_0,\sigma\right) \simeq \Ext_H^i\left(\ind_Q^H\left(\ind_{U_{k(a-1)}}^R\psi_{k(a-1)}\right)\boxtimes \pi_0,\sigma\right).
\] 
On the other hand, by the induction hypothesis $\Ext^i_H\left(\xi\rtimes\pi_0,\sigma\right) = 0$ for all $i\geq 1$. Then (\ref{E:AuxVanish3}) follows from the combination of all these isomorphisms. Applying the functor $\Hom_H(-,\sigma)$ to (\ref{E:S.E.S-3}) we get
\[
    \Ext^i_H\left(\pi_1\oplus\pi_2,\sigma\right) \simeq  \Ext^{i+1}_H\left(\pi',\sigma\right) 
\]
for all $i \geq 1$.

\vskip 5pt

Dually, applying both the contragredient and MVW-involution to (\ref{E:S.E.S-3}), we get a dualized short exact sequence
\begin{equation*}\label{E:S.E.S-4}\tag{$\diamondsuit$}
    0 \lra \pi' \lra \delta_\rho(x-1, -x)\rtimes \pi_0 \lra \pi_1\oplus\pi_2 \lra 0,
\end{equation*}
Using the same Mackey argument as we have done several times, we claim that
\begin{equation*}\label{E:AuxVanish4}\tag{$\ddagger.2$}
    \Ext_H^i\left(\delta_\rho(x-1, -x)\rtimes \pi_0, \sigma\right) = 0
\end{equation*}
for all $i\geq 1$. Indeed, the $H$-subrepresentation $\Pi'_{op}$ of $\delta_\rho(x-1, -x)\rtimes \pi_0$ corresponding to the open orbit can be computed totally the same as $\Pi_{op}$; as for the quotient $\Pi'_c$ corresponding to the closed orbit, we have 
\[
    \Pi'_c = \delta_\rho\left(x-\frac{1}{2}, -x+\frac{1}{2}\right)\rtimes \left(\pi_0~\big|_{H_0}\right).
\]
To see $\Pi'_c$ will not contribute to Ext spaces, again we appeal to the second Frobenius reciprocity and K\"unneth's formula. We have 
\[
    \Ext^i_H\left(\Pi'_c, \sigma\right) = \Ext^i_{\GL_{k(a-1)}\times H_0}\left(\delta_\rho\left(x-\frac{1}{2}, -x+\frac{1}{2}\right)\boxtimes \left(\pi_0~\big|_{H_0}\right), \Jac_{\overline{Q_c}}\sigma\right),
\]
where $Q_c$ is the stabilizer of a representative in the closed orbit, which is also a maximal parabolic subgroup of $H$ with Levi component $\GL_{k(a-1)}\times H_0$. Recall that $\sigma = \lambda\rtimes \sigma_{bp}$. By the Tadi\'c's formula, Lemma \ref{L:CasselmanCri} and assumption (\textbf{A.3}), if we write
\[
    s.s.\Jac_{\overline{Q_c}}\sigma = \sum_{j} \tau_j\boxtimes \sigma_j
\]
for some irreducible representations $\tau_j$ and $\sigma_j$ of $\GL_{k(a-1)}$ and $H_0$, then for each $\tau_j$, either it has strictly negative exponent, or it does not live in the same cuspidal component as $\delta_\rho\left(x-\frac{1}{2}, -x+\frac{1}{2}\right)$. In conclusion we always have
\[
    \Ext^i_{\GL_{k(a-1)}}\left(\delta_\rho\left(x-\frac{1}{2}, -x+\frac{1}{2}\right) , \tau_j\right) = 0
\]
for all $i\geq 0$. Hence we deduce from the K\"unneth's formula that
\[
    \Ext^i_{H}\left(\Pi'_c, \sigma\right) = 0
\]
for all $i\geq 0$. Applying the functor $\Hom_H(-,\sigma)$ to the short exact sequence (\ref{E:S.E.S-4}), we deduce that
\[
    \Ext^i_H\left(\pi',\sigma\right) \simeq \Ext^{i+1}_H\left(\pi_1\oplus\pi_2,\sigma\right)
\]
for all $i \geq 1$. Thus one can see that $\Ext^i_H(\pi_1\oplus\pi_2,\sigma)$ is periodic, hence must vanishes.

\vskip 10pt

So now we have completed the proof of Proposition \ref{P:discCase}.

\vskip 10pt

\section{Proof of the main theorem II: tempered case}

In this section we will finish the proof of Theorem \ref{T:Main} based on Proposition \ref{P:discCase}. Let $M$ be a tempered L-parameter of $G$ and $N$ of $H$. If $M$ does not satisfy (\textbf{A.0}), then there is an irreducible $k$-dimensional representation $\rho$ of $W_E$ and $a\in\Z_{>0}$, such that
\[
    M = \rho S_a + M_0 + \left(\rho^c S_a\right)^\vee.
\]
Here $M_0$ is a tempered L-parameter for $G_0=G(V_0)$, with $V_0\subset V$ a non degenerate subspace such that $V\simeq V_0+\calH^{ka}$. Then by the LLC we know that for any $\pi\in\Pi_M(G)$, there exists a unique $\pi_0\in\Pi_{M_0}(G_0)$ such that
\[
    \pi\hookrightarrow \delta_\rho\left(\frac{a-1}{2}, -\frac{a-1}{2}\right)\rtimes \pi_0
\]
as a direct summand. Hence to show that $\Ext^i_H\left(\pi,\sigma\right) = 0$ for all $\sigma\in\Pi_N(H)$ and $i\geq 1$, it suffices to show that
\[
    \Ext^i_H\left(\delta_\rho\left(\frac{a-1}{2}, -\frac{a-1}{2}\right)\rtimes \pi_0,\sigma\right) = 0.
\]
Same as the argument in Section \ref{SS:Induction}, one can choose a sequence of non conjugate self-dual characters $\chi_1,\cdots,\chi_{ka}$ suitably, such that
\[
    \Ext^i_H\left(\delta_\rho\left(\frac{a-1}{2}, -\frac{a-1}{2}\right)\rtimes \pi_0,\sigma\right) \simeq \Ext^i_H\left(\chi_1\times\cdots\times\chi_{ka}\rtimes \pi_0,\sigma\right).
\]
Note that $\chi_1\times\cdots\times\chi_{ka}\rtimes \pi_0$ has L-parameter
\[
    \left(\chi_1+\cdots+\chi_{ka}\right) + M_0 + \left(\chi_1^c+\cdots+\chi_{ka}^c\right)^\vee.
\]
Repeating this process, after finite steps we can replace $M$ by another L-parameter $M^\flat$ of $G$ which satisfies (\textbf{A.0}). Similarly, using the same argument and Corollary \ref{C:Switch}, we can also replace $N$ by another L-parameter $N^\flat$ of $H$ which satisfies (\textbf{A.0}). Then the desired conclusion follows from Proposition \ref{P:discCase}.

\vskip 5pt

This completes the proof of Theorem \ref{T:Main}.

\vskip 5pt

\begin{Rmk}\label{R:Highcodim}
Up to now we have only considered the case when $W\subset V$ is a codimension $1$ subspace. This is the so called ``basic Bessel case'' of the Gan--Gross--Prasad model. In general, if $W\subset V$ is a non degenerate subspace of odd codimension such that $V/W$ is split, then following \cite[Sect. 12]{MR3202556} one can define the Bessel subgroup 
\[
    {\rm Bes}(W,V) \subset H(W)\times G(V)
\]
and a distinguished character $\mu$ of ${\rm Bes}(W,V)$, and consider Ext spaces
\[
    \Ext^i_{{\rm Bes}(W,V)}\left(\pi,\sigma\right):=\Ext^i_{{\rm Bes}(W,V)}\left(\pi\boxtimes \sigma^\vee,\mu\right)
\]
for representations $\pi$ of $G(V)$ and $\sigma$ of $H(W)$. Using the same argument as \cite[Thm. 15.1]{MR3202556}, one can show that Theorem \ref{T:Main} implies the same Ext-vanishing result for higher codimension cases, namely, if $\pi$ and $\sigma$ are tempered, then
\[
    \Ext^i_{{\rm Bes}(W,V)}\left(\pi,\sigma\right)=0
\]
for all $i\geq 1$. There are also other cases called ``Fourier--Jacobi case'' and ``twisted Fourier--Jacobi case''. We will investigate them in the future.
\end{Rmk}

\vskip 10pt

\section{Proof of the projectivity result}

In this section we prove the byproduct: Theorem \ref{T:ProjSteinberg}. So let $\pi$ be an irreducible tempered representation of $G$, which is of maximal derivative as in Definition \ref{D:GSteinberg}. Assume that $\pi$ is non cuspidal since otherwise the desired conclusion obviously holds. We fix an irreducible unitary supercuspidal representation $\rho_0$ of $\GL_{k_0}$ and $x_0$ a real number such that 
\[
    \Jac_{\rho_0,x_0}\pi \neq 0.
\]
Then $\rho_0$ and $2x_0+1$ satisfy the requirement in Definition \ref{D:GSteinberg}, and we have $x\geq \frac{\dim V-2k_0}{2k_0}$.

\vskip 5pt

According to the result of Aizenbud-Sayag \cite[Thm. 3.2]{MR4054816}, the restriction of $\pi$ to $H$ is a so called ``locally finite representation'' (see also \cite[Thm. 1.2]{prasad2023homological}). Then, by a general result of Chan--Savin \cite[Thm. A.1]{MR3910471}, to show that $\pi~\big|_H$ is projective, it remains to show that
\begin{equation*}\label{E:ProjSteinbergExt}\tag{$\natural$}
    \Ext_H^i\left(\pi,\sigma\right) = 0
\end{equation*}
for all irreducible representation $\sigma$ of $H$ and all $i\geq 1$. We show this by estimating the exponents of Jacquet modules of $\sigma$. 

\vskip 5pt

Apparently, there are two possibilities of $\sigma$: either it is tempered, or it is non-tempered. In the tempered case, by \cite[Cor. 2.9, Lem. 4.1]{MR4055180} the exponents of Jacquet modules of $\sigma$ are bounded by the L-parameter of $\sigma$: if there is some irreducible unitary supercuspidal representation $\rho$ of $\GL_k$ such that $\Jac_{\rho,x}\sigma\neq 0$ for some real number $x$, then $x\in\frac{1}{2}\Z_{\geq 0}$ and
\[
    \rho S_{2x+1}\subset LL\left(\sigma\right),
\]
which in particular implies that $x\leq \frac{d_H-k}{2k}$. Here we recall that $d_H$ is the dimension of L-parameters of $H$. In the non-tempered case, by the Langlands classification one can write $\sigma$ as the unique irreducible subrepresentation of a parabolic induction of the form
\[
    \delta_{\rho_1}(x_1,y_1)\times \delta_{\rho_2}(x_2,y_2)\times\cdots\times \delta_{\rho_r}(x_r,y_r)\rtimes \sigma_0,
\]
where each $\rho_i$ is an irreducible unitary representation of a general linear group $\GL_{k_i}$, $x_i,y_i$ are real numbers such that
\[
    x_1+y_1\leq x_2+y_2\leq\cdots\leq x_r+y_r< 0,
\]
and $\sigma_0$ is an irreducible tempered representation of $H_0=H\left(W_0\right)$, with $W_0\subset W$ a non degenerate subspace of proper dimension. Note that
\[
    \sum_{i=1}^r k_i(x_i-y_i+1)\leq r_W,
\]
where $r_W$ is the Witt index of $W$. In particular, we deduce that $x_1< \frac{r_W-k_1}{2k_1}$, and also $\Jac_{\rho_1,x_1}\sigma\neq 0$. The same analysis is also applicable to the Zelevinsky--Aubert dual $\widehat{\sigma}$ of $\sigma$. Then by \cite[Prop. 3.9(1)]{MR4549708}, we get the following trivial lower bound of the largest exponent of Jacquet modules:

\begin{Lem}
In the context of above discussion, as long as $\sigma$ is not supercuspidal, there exists an irreducible unitary supercuspidal representation $\rho$ of a general linear group $\GL_k$, such that $\Jac_{\rho,x}\sigma\neq 0$ for some real number 
    \[
        x \geq \frac{k-d_H}{2k}.
    \]
\end{Lem}

\vskip 5pt

With this estimation at hand, now we prove (\ref{E:ProjSteinbergExt}). If $\sigma$ is supercuspidal, then (\ref{E:ProjSteinbergExt}) holds since $\sigma$ is compact. So from now on we assume that $\sigma$ is not supercuspidal. Then by the lemma above, we can embed $\sigma$ into an induced representation
\[
    \sigma\hookrightarrow\rho\nu^x\rtimes \sigma_1,
\] 
where $\rho$ is an irreducible unitary supercuspidal representation of a general linear group $\GL_k$, $x\geq \frac{k-d_H}{2k}$ is a real number, and
$\sigma_1$ is an irreducible representation of $H_1=H\left(W_1\right)$, with $W_1\subset W$ a non degenerate subspace such that $W\simeq W_1+\calH^{k}$. Let $\chi$ be a unitary non conjugate self-dual character of $\GL_1$ which is not an unramified twist of $\rho^c$, and does not occur in the cuspidal component of $\sigma_1^\vee$. Applying the MVW-involution and the parabolic induction functor $\chi\rtimes-$, we obtain
\[
    \chi\rtimes\sigma^\vee\hookrightarrow \rho^c\nu^x\rtimes \left(\chi\rtimes\sigma_1^\vee\right).
\]
Then by Corollary \ref{C:Switch} and a standard degree shifting argument, in order to show
\[
    \Ext_H^i\left(\pi,\sigma\right)\simeq \Ext_G^i\left(\chi\rtimes\sigma^\vee,\pi^\vee\right)
\]
are zero for all $i\geq 1$, it suffices to show that
\[
    \Ext_G^i\left(\rho^c\nu^x\rtimes \left(\chi\rtimes\sigma_1^\vee\right),\pi^\vee\right) = 0
\]
for all $i\geq 1$. Once again we appeal to the Mackey theory. The quotient of $\rho^c\nu^x\rtimes \left(\chi\rtimes\sigma_1^\vee\right)$ corresponding to the closed orbit is
\[
    \rho^c\nu^{x+\frac{1}{2}}\rtimes \left(\chi\rtimes\sigma_1^\vee~\big|_{G_1}\right),
\]
where $G_1=G\left(V_1\right)$, with $V_1\subset V$ a non degenerate subspace such that $V\simeq V_1+\calH^k$. We compute
\[
    \Ext_G^i\left(\rho^c\nu^{x+\frac{1}{2}}\rtimes \left(\chi\rtimes\sigma_1^\vee~\big|_{G_1}\right),\pi^\vee\right)\simeq \Ext_G^i\left(\pi,\left(\rho^c\right)^\vee\nu^{-x-\frac{1}{2}}\rtimes \left(\chi\rtimes\sigma_1^\vee~\big|_{G_1}\right)^\vee\right).
\]
If $\rho\notin\left\{\rho_0, \left(\rho_0^c\right)^\vee\right\}$, then clearly the closed orbit does not contribute to Ext spaces. If $\rho\in\left\{\rho_0, \left(\rho_0^c\right)^\vee\right\}$ (so in this case $k=k_0$), we have
\[
    -x-\frac{1}{2}\leq \frac{d_H-k}{2k}-\frac{1}{2}< \frac{\dim V-2k_0}{2k_0}.
\] 
It follows from our assumption on $\pi$ that the contribution of the closed orbit is also zero. %We claim that the equality can not hold. Indeed, if the equality holds, then $d_G=d_H$ and $d_H= k\left(-2x+1\right)$.
Therefore only the open orbit has contribution. Let $\chi_1,\cdots,\chi_k$ be a sequence of unitary non conjugate self-dual characters of $\GL_1$, such that $\chi_j$ does not occur in the cuspidal components of $\pi^\vee$ and $\sigma_1^\vee$ for all $j\in\{1,\cdots,k\}$. Similar to the argument in Section \ref{SS:Induction}, we have
\[
    \Ext_G^i\left(\rho^c\nu^x\rtimes \left(\chi\rtimes\sigma_1^\vee\right),\pi^\vee\right) \simeq \Ext_G^i\left(\chi_1\times\cdots\times\chi_k\rtimes \left(\chi\rtimes\sigma_1^\vee\right),\pi^\vee\right).
\]
Now, if $\sigma_1$ is supercuspidal, then $\chi_1\times\cdots\times\chi_k\rtimes \left(\chi\rtimes\sigma_1^\vee\right)$ is tempered, which implies that 
\[
    \Ext_G^i\left(\chi_1\times\cdots\times\chi_k\rtimes \left(\chi\rtimes\sigma_1^\vee\right),\pi^\vee\right) = 0
\]
for all $i\geq 1$ by Theorem \ref{T:Main}. If $\sigma_1$ is not supercuspidal, we can repeat this process to the representation $\sigma_1$. After finite steps we can eventually replace $\chi\rtimes \sigma^\vee$ by a tempered representation of the form
\[
    \chi_1\times\cdots\times\chi_k\times\cdots\times\chi_s\rtimes\left(\chi\rtimes\sigma_{sc}^\vee\right),
\]
where for each $j$, $\chi_j$ is a unitary character of $\GL_1$, and $\sigma_{sc}$ is a supercuspidal representation. Then the desired Ext-vanishing result (\ref{E:ProjSteinbergExt}) follows from our main result Theorem \ref{T:Main}.

\vskip 5pt

This completes the proof of Theorem \ref{T:ProjSteinberg}.

\vskip 10pt

\section*{Acknowledgments}

The author would like to thank Prof. Kei Yuen Chan, Prof. Dipendra Prasad, Dr. Chuijia Wang and Dr. Jialiang Zou for helpful discussions.

\vskip 15pt

\bibliographystyle{alpha}
%\nocite{*}
\bibliography{ExtBVanish4TRef}

\end{document}